\documentclass[12pt]{amsart}
\usepackage{verbatim}
\usepackage{amsmath}
\usepackage{amssymb}

\numberwithin{equation}{section}
 
\newtheorem{theorem}[equation]{Theorem}
\newtheorem{lemma}[equation]{Lemma}

\newtheorem{proposition}[equation]{Proposition}
\newtheorem{corollary}[equation]{Corollary}

\newtheorem*{theorem*}{Theorem}
 
\theoremstyle{definition}
\newtheorem*{definition}{Definition}

\theoremstyle{remark}

\newcommand{\m}{{\mathfrak{m}}}

\newcommand{\U}{{\mathcal U}}

\newcommand{\A}{{\mathcal A}}
\newcommand{\B}{{\mathcal B}}
\newcommand{\C}{{\mathcal C}}

\newcommand{\F}{{\mathcal F}}

\newcommand{\Xhat}{{\hat{X}}}
\newcommand{\Chat}{{\hat{\C}}}

\newcommand{\onto}{\twoheadrightarrow}  
\newcommand{\into}{\hookrightarrow}     

\def\CC{{\mathbb C}}

\def\RR{{\mathbb R}}
\def\ZZ{{\mathbb Z}}
\def\PP{{\mathbb P}}

\newcommand{\Max}{\operatorname{Max}}

\newcommand{\Kstar}{K^*}

\newcommand{\Star}{\operatorname{B}}

\newcommand{\an}{{\rm an}}

\begin{document}

\title{Uniform Structures and Berkovich Spaces}


\author{Matthew Baker}

\begin{abstract}
A uniform space is a topological space together with some additional structure which allows one to 
make sense of uniform properties such as completeness or uniform convergence.
Motivated by previous work of J. Rivera-Letelier, we give a new construction of the Berkovich analytic space associated to an affinoid algebra 
as the completion of a canonical uniform structure on the associated rigid-analytic space.
\end{abstract}

\thanks{I would like to thank Xander Faber for helpful discussions, and for 
his assistance with the proof of Proposition~\ref{ShrinkingLaurentProp}.
I would also like to thank Juan Rivera-Letelier for telling me about the theory 
of uniform spaces, and for sharing his preprint \cite{RLUniform}.  
The exposition in this paper was significantly improved by the comments of
Brian Conrad, Rob Benedetto, Xander Faber, and Juan Rivera-Letelier.
This work was supported by NSF grants DMS-0300784 and DMS-0600027.}
                                       
\maketitle


\section{Introduction}


In \cite{RLUniform}, Juan Rivera-Letelier proves that there is a canonical uniform structure on the projective
line over an arbitrary complete non-archimedean ground field $K$ whose completion is
the Berkovich projective line over $K$.  In this note, we extend Rivera-Letelier's ideas by proving
that there is a canonical uniform structure on any rigid-analytic affinoid space $X$, defined in a natural way 
in terms of the underlying affinoid algebra $\A$, whose completion yields the analytic space
$X^{\an}$ associated to $X$ by Berkovich in \cite{BerkovichBook}.
Actually, the existence of such a uniform structure on $X$ can be deduced easily from
general properties of Berkovich spaces and uniform spaces
(see \S\ref{BerkovichSection}), so the interesting point here is
that one can define the canonical uniform structure on $X$ and verify its salient properties
using only some basic facts from rigid analysis, and without invoking Berkovich's theory.
Therefore one arrives at a new method for constructing (as a topological space) the 
Berkovich space associated to $X$, and for 
verifying some of its important properties (e.g., compactness).

\medskip

A uniform space is a set endowed with a ``uniform structure''.  
Uniform spaces are intermediate objects between topological spaces and metric spaces; 
in particular, every uniform space has a natural topology, and every metric space 
is endowed with a canonical uniform structure.  
Every topological group also has a canonical uniform structure, so uniform spaces are 
more general than metric spaces.

Intuitively, in a uniform space one can formalize the assertion ``$x$ is as close to $y$ as $x'$ is
to $y'$'', while in a general topological space one can only make sense of ``$x$ is as close to $y$ as 
$x'$ is to $y$''.
Thus one can formalize uniform properties such as completeness, uniform continuity, and uniform convergence in a uniform space but not in a
general topological space.  

\medskip

As an illustrative example, and one which is relevant to the larger context of this paper, 
consider the topological space $\ZZ_p$ of $p$-adic integers.
The space $\ZZ_p$ is metrizable, of course, but there are many distinct metrics which
give rise to the same $p$-adic topology.  On the other hand, there is one and only one uniform structure on $\ZZ_p$ which gives rise to the $p$-adic topology,
and restricting this to $\ZZ$ yields a canonical uniform structure on the set of integers whose completion
yields the $p$-adic integers.

\medskip

We now explain how Berkovich's theory of non-archimedean analytic spaces can be viewed in an analogous way.
In rigid analytic geometry, one works with analogues of manifolds 
(or more generally of real and complex analytic spaces)
where the field $\CC$ of complex numbers is replaced by a complete non-archimedean field $K$. 
To deal with the fact that the canonical topology on an analytic manifold $X$ over $K$ is
totally disconnected and (when $K$ is not a local field) not even locally compact,
John Tate introduced rigid analytic spaces, whose building blocks are
affinoid spaces of the form $\Max(\A)$,
where $\A$ is an affinoid algebra over $K$ and $\Max(\A)$ is the space of maximal ideals in $\A$.
Tate dealt with the problems which arise in the non-archimedean world by introducing a Grothendieck topology
on $\Max(\A)$.
A number of years later, Vladimir Berkovich had the insight that 
one could instead look at the larger topological space $M(\A)$ consisting of all bounded multiplicative seminorms on $\A$, together with a
natural topology.  There is a continuous embedding $\Max(\A) \into M(\A)$ having dense image, but the space $M(\A)$ has many extra
points which allow one to recover the same sheaf theory as in Tate's theory, but using actual
topological spaces, as opposed to just a Grothendieck topology. One can glue the affinoid spaces $M(\A)$ together
to obtain a global analytic space $X^{\an}$ functorially associated to any ``reasonable'' rigid analytic space $X/K$.
In particular, every algebraic variety $X/K$ has a Berkovich analytification $X^{\an}$.
The analytic spaces in Berkovich's theory are very nice from a topological point of view.
For example, Berkovich affinoids and the analytifications of projective varieties over $K$ are compact Hausdorff spaces.
In addition, a Berkovich affinoid is path-connected if and only if the corresponding rigid space is connected (in the sense
of rigid analysis).


\medskip

In this paper, we will show that given an affinoid space $X = \Max(\A)$ in the sense of Tate, 
one can construct the corresponding Berkovich space $X^{\an}$ as the completion of a canonical
uniform structure on $X$ generated by finite coverings of $X$ by {\em open Laurent domains}.
The fact that this really defines a uniform structure can be deduced from Berkovich's theory, but it can also be proved directly 
using rigid analysis, specifically the Maximum Modulus Principle for affinoid algebras.
This allows one to prove, independently of Berkovich's theory and without ever mentioning multiplicative seminorms,
that $X$ can be embedded in a natural way as a dense subspace of a compact Hausdorff topological space $\hat{X}$.
General considerations then allow one to conclude {\em a posteriori} that $\hat{X}$ and $X^{\an}$ are homeomorphic.



\medskip

{\bf Remarks:}

\medskip

1. For the purposes of this paper, Berkovich spaces are considered only as topological spaces.  Our reason for
ignoring sheaves and rings of functions on these spaces is simply that we have nothing new to add to what 
is already done by Berkovich in \cite{BerkovichBook} and \cite{BerkovichIHES}.
Also, we should remark that the affinoid spaces we consider in this paper are what Berkovich calls {\em strict affinoids}. 

\medskip

2. The uniform structure on Berkovich spaces is also interesting for certain applications.
For example, there is a rich theory of
rational iteration on the Berkovich projective line (see, e.g., \cite{RLFatouJulia})
which is analogous to the classical Fatou-Julia theory on $\PP^1(\CC)$.
The Fatou set is defined classically as the locus of equicontinuity or normality of 
the family of iterates of $\varphi$; both of these notions make sense in a uniform space but neither makes sense in a general topological space.
So in order to give a definition of the Berkovich Fatou and Julia sets of a rational map $\varphi$ which mimics the classical definitions, one 
should use the uniform structure on $(\PP^1)^{\an}$.  

\medskip

3. It may be possible to gain insight into the structure of higher-dimensional Berkovich spaces by using the fact that
points in the completion of a uniform space correspond to minimal Cauchy filters, and thus
points of a Berkovich space can be thought of as minimal Cauchy filters on the underlying rigid-analytic space.  
For example, when $K$ is algebraically closed, 
Rivera-Letelier uses minimal Cauchy filters to give a new derivation of Berkovich's classification 
of points of $(\PP^1)^{\an}$ in terms of nested sequences of closed disks in $\PP^1(K)$
(see \cite{RLUniform}).

\subsection{Organization of the paper}

This paper is organized as follows.  
In \S\ref{UniformSpaceSection}, we will give an overview of the theory of uniform spaces, emphasizing the definition in
terms of uniform covers.  In \S\ref{StarRefinementSection}, we formulate and prove a technical result (the ``Shrinking Lemma'') 
which simplifies the task of proving that a given collection of coverings satisfies the axioms of a uniform structure.  
After recalling the necessary background from rigid analysis, we prove 
in \S\ref{RigidAnalysisSection} that finite coverings by open Laurent domains 
satisfy the conditions of the shrinking lemma.  This allows us to prove directly, 
without utilizing Berkovich's theory, that every rigid-analytic affinoid $X$ admits a canonical uniform structure.
After recalling some basic facts about Berkovich spaces, we then prove that the completion of this uniform structure
coincides with the Berkovich analytic space associated to $X$.  Finally, in a brief appendix, we describe the
``entourage'' approach to uniform spaces, and indicate how one translates between the uniform cover and entourage
points of view.

\section{Uniform spaces}
\label{UniformSpaceSection}

Uniform spaces were first introduced by Andr{\'e} Weil, and there are several equivalent definitions in use today. 
For example, Bourbaki \cite{BourbakiGT} defines uniform structures in terms of ``entourages'', while John Tukey
\cite{Tukey} defines them in terms of ``uniform covers''.  There is also a definition due to Weil in terms of 
``pseudometrics''.  We will utilize the uniform cover definition of Tukey throughout this paper.
A general reference for this section is Norman Howes' book \cite{Howes}.
A more detailed overview of the theory of uniform spaces can be found in 
the online encyclopedia {\em Wikipedia} (\verb+http://en.wikipedia.org/wiki/Uniform_space+), 
and {\em The Encyclopaedia of Mathematics}, available online at 
\linebreak
\verb+http://eom.springer.de/u/u095250+.  For a dictionary between the uniform cover and 
entourage approaches, see the Appendix.

\subsection{Uniform covers}
\label{UniformCoverSection}

The definition of uniform spaces in terms of {\em uniform covers} is as follows.

Let $X$ be a set.  A {\em covering} (or cover) of $X$ is a collection of subsets of $X$ whose union is $X$.
A covering $C$ {\em refines} a covering $D$, written $C < D$, if every $U \in C$ is contained in some $V \in D$.
A covering $C$ {\em star-refines} a covering $D$
if for every $U \in C$, the union of all elements of $C$ meeting $U$ is contained in some element of $D$.
For example, if $X$ is a metric space and $\epsilon > 0$, the collection $C$ of all open balls of radius $\epsilon/3$ in
$X$ star-refines the collection $D$ of all open balls of radius $\epsilon$.

The {\em intersection} $C \cap D$ of two coverings is the covering $\{ U \cap V \; | \; U \in C, \, V \in D \}$.

A {\em uniform structure} on a set $X$ is a collection $\C$ of coverings of $X$, called the {\em uniform coverings}, 
satisfying the following two axioms:

\begin{itemize}
\item[(C1)] If $C \in \C$ and $C < D$, then $D \in \C$. 
\item[(C2)] Given $D_1,D_2 \in \C$, there exists $C \in \C$ which star-refines both
$D_1$ and $D_2$.
\end{itemize}

It is easy to verify that in the presence of axiom (C1), axiom (C2) is equivalent to the following
two axioms:

\begin{itemize}
\item[(C2)(a)] If $C_1,C_2 \in \C$, then $C_1 \cap C_2 \in \C$.
\item[(C2)(b)] Every $C \in \C$ has a star-refinement in $\C$.
\end{itemize}

If a collection $\C'$ of coverings of $X$ satisfies (C2) but not (C1), we call it a {\em base} for a uniform structure.
If $\C'$ is a base for a uniform structure, the {\em uniform structure generated by} $\C'$ is
the set of all coverings of $X$ having a refinement in $\C'$.

A {\em uniform space} is a pair $(X,\C)$, where $\C$ is a uniform structure on the set $X$.
A uniform space is called {\em Hausdorff} if it satisfies the following additional axiom:
\begin{itemize}
\item[(C3)] For any pair $x,y$ of distinct points in $X$, there exists $C \in \C$ such that 
no element of $C$ contains both $x$ and $y$.
\end{itemize}

If $X$ is a metric space, one can define a canonical Hausdorff uniform structure on $X$ by taking as a base for the uniform covers
of $X$ the set $\C' = \{ C_\epsilon \; | \; \epsilon > 0 \}$, where 
$C_\epsilon$ consists of all open balls of
radius $\epsilon$ in $X$.  The fact that $\C'$ satisfies axiom (C2) follows from the triangle inequality, 
which guarantees that $C_{\epsilon / 3}$ is a star-refinement of $C_{\epsilon}$.
Axiom (C2) can thus be thought of as a version of the triangle inequality which makes sense without the presence
of a metric.

If $(X,\C)$ is a uniform space, then for $x \in X$ and $C \in \C$, one
defines $B(x,C)$, the ``ball of size $C$ around $x$'', to be the union of all elements of $C$ containing $x$.
By definition, the {\em uniform topology} on $(X,\C)$ is the topology for which 
$\{ B(x,C) \; | \; C \in \C \}$ forms a fundamental system of neighborhoods of $x \in X$.
In other words, a subset $U \subseteq X$ is open in the uniform topology if and only if
for each $x \in U$, there exists a cover $C \in \C$ with $B(x,C) \subseteq U$.
If $X$ is a metric space, then $B(x,C_{\epsilon/2})$ contains each open ball of radius $\epsilon/2$ about $x$ and is contained in 
the open ball of radius $\epsilon$ about $x$, and thus the 
metric topology on $X$ coincides with the uniform topology on $X$.
The uniform topology on $(X,\C)$ is Hausdorff if and only $(X,\C)$ satisfies axiom (C3).

A mapping $f : (X,\C_X) \to (Y,\C_Y)$ between uniform spaces is called {\em uniformly continuous} if for
any $C \in \C_Y$, the covering $\{ f^{-1}(U) \; | \; U \in C \}$ 
belongs to $\C_X$.
A uniformly continuous mapping is continuous relative to the uniform topologies on $X$ and $Y$, and
if $X$ and $Y$ are metric spaces (endowed with the canonical uniform structure), 
then $f$ is uniformly continuous in the sense of uniform spaces if and only if it
is uniformly continuous in the sense of metric spaces.
Uniform spaces, together with the uniformly continuous mappings between them, form a category.

A uniform space $X$ is called a {\em uniform subspace} of the uniform space $Y$ if 
there is a uniformly continuous injective map $i : X \into Y$.  There is a unique uniform
structure on any subset $A$ of $Y$ making $A$ into a uniform subspace of $Y$.
We say that $X$ is a {\em dense} uniform subspace of $Y$ if it is dense with respect to 
the uniform topology.

\subsection{Completeness and compactness}
\label{CompletenessSection}

In this section, we discuss the uniform notion of completeness and its relation to compactness.

First, we recall the definition of a filter. A {\em filter} on a set $X$ is a collection $\F$ of subsets of $X$
such that:
\begin{itemize}
\item[(F1)] $\emptyset \not\in \F$. 
\item[(F2)] If $U_1,U_2 \in \F$ then $U_1 \cap U_2 \in \F$.
\item[(F3)] If $U \in \F$ and $U \subseteq V$, then $V \in \F$.
\end{itemize}
The collection of filters on $X$ is partially ordered by inclusion: we write $\F \leq \F'$ if
every $U \in \F$ is contained in some $U' \in \F'$.  In this case, $\F$ is said to {\em refine} $\F'$.

As an example, if $X$ is a topological space and $x \in X$, the set of all neighborhoods of $x$ forms a filter $N_x$, 
called the {\em neighborhood filter} of $x$.
We say that a filter $\F$ {\em converges} to $x$ (or $x$ is a {\em limit point} of $\F$) if $N_x \leq \F$.
It is easy to show that if $X$ is Hausdorff, then a filter on $X$ can have at most one limit point.

As another example, if $(x_n)$ is a sequence of elements of $X$, the collection of all subsets $V$ of $X$ such that
$x_n \in V$ for all sufficiently large $n$ is a filter.  In this way, filters generalize sequences.

A filter on a uniform space $(X,\C)$ is called {\em Cauchy} if it contains some element of each uniform covering.
A uniform space $(X,\C)$ is called {\em complete} if every Cauchy filter converges to a point of $X$.
This definition coincides with the metric definition of completeness when $X$ is a metric space.

If $(X,\C)$ is a Hausdorff uniform space, there is a 
complete Hausdorff uniform space $(\Xhat,\Chat)$, unique up to uniform isomorphism,
which contains $X$ as a dense uniform subspace.  The space
$(\Xhat,\Chat)$ is called the {\em completion} of $(X,\C)$.


Generalizing the usual construction of the completion of a metric space,
the completion of a Hausdorff uniform space $(X,\C)$ can be defined as the set of 
equivalence classes of Cauchy filters, endowed with a natural topology.
Here two Cauchy filters are called {\em equivalent} if their intersection is also a Cauchy filter.

Alternatively, one can show that each equivalence class of Cauchy filters contains a
unique minimal element with respect to inclusion, so that $\Xhat$ can be defined on the level of points as the set of 
{\em minimal Cauchy filters} on $X$.
For example, the neighborhood filter $N_x$ of a point $x \in X$ is a minimal Cauchy filter.  
This provides a concrete description of the natural embedding $X \into \Xhat$.

It is important to note that a set $X$ might have many different uniform structures which all induce
the same topology, and the completions of these uniform structures will in general be different.
For example, if $K$ is a complete non-archimedean field, then the closed unit disc $D$ in $K$ is complete with respect to the uniform structure defined by the 
usual metric.  On the other hand, as we will see, there is another uniform structure on $D$, defined in terms of finite coverings by
open Laurent domains, which gives rise to the same topology, but whose completion is the Berkovich analytic space 
associated to $D$.

A uniform space $(X,\C)$ is called {\em totally bounded} (or {\em precompact}) if each uniform cover has a finite
subcover.  A basic fact about Hausdorff uniform spaces is that $(X,\C)$ is compact in its uniform topology if and only if
$(X,\C)$ is complete and totally bounded.  
As a consequence, one sees that the completion of a Hausdorff uniform
space $(X,\C)$ is compact if and only if $(X,\C)$ is totally bounded.

Another basic result in the theory of uniform spaces is that if $X$ is a compact Hausdorff topological space, 
then there is a {\em unique} uniform structure on $X$ compatible with the given topology.  A base for the 
uniform structure on $X$ is given by the collection of all coverings of $X$ by open subsets.
It is not hard to see that if $\B$ is any base for the topology on $X$, 
a base for the uniform structure on $X$ is also given by the collection $\C_\B$ of all finite coverings of $X$
by elements of $\B$.

\section{Existence of star-refinements}
\label{StarRefinementSection}

It can be difficult in practice to verify axiom (C2) (or axiom (C2)(b)).
In this section, we present a criterion (the ``Shrinking Lemma'') for the existence of star-refinements which
allows one to verify axiom (C2)(b) rather painlessly.  This criterion will be used in the next
section to prove that every rigid-analytic affinoid space has a canonical uniform structure.

First, we need some preliminary facts and definitions.

Let $X$ be a set, and let $C$ be a covering of $X$.  
If $H \subseteq X$, we define $\Star(H,C)$ to be the union of all elements of $C$
meeting $H$.
For example, if $x \in X$ then $\Star(x,C)$ 
is just the ball of size $C$ around $x$
in the terminology of \S\ref{UniformCoverSection}.
We define $C^* = \{ \Star(U,C) \; | \; U \in C \}$ and 
$C^\Delta = \{ \Star(x,C) \; | \; x \in X \}$.
It is easy to check that $C^*$ and $C^\Delta$ are both coverings of $X$.
A straightforward argument shows that
\begin{equation}
\label{eq:StarDelta}
C^\Delta < C^* < \left(C^{\Delta}\right)^{\Delta} \ .
\end{equation}
It is also easy to see that if $C < D$, then $C^* < D^*$ and 
$C^\Delta < D^{\Delta}$.
(Tukey \cite{Tukey} calls these various assertions about stars and deltas
the ``calculus of coverings''.)

By definition, one sees that a covering $C$ star-refines a covering $D$ if and only if $C^* < D$.
Similarly, we will say that $C$ {\em delta-refines} $D$ if $C^\Delta < D$.
As an example, if $X$ is a metric space and $\epsilon > 0$, the collection $C$ of all balls of radius $\epsilon/2$ in
$X$ delta-refines the collection $D$ of all balls of radius $\epsilon$.

By (\ref{eq:StarDelta}), it is easy to see that axiom (C2) in the definition of a uniform structure is equivalent to:
\begin{itemize}
\item[$({\rm C2})^{\prime}$] Given $D_1,D_2 \in \C$, there exists $C \in \C$ which delta-refines both $D_1$ and $D_2$.
\end{itemize}
Similarly, axiom (C2)(b) is equivalent to the {\it a priori} weaker condition
\begin{itemize}
\item[$({\rm C2})^{\prime}({\rm b})$] Every cover in $\C$ has a delta-refinement in $\C$.
\end{itemize}

A covering $C$ is called {\em finite} if it has finitely many elements.
If $C,D$ are finite coverings of $X$ with $C<D$, we say (following \cite[p. 57]{Howes}) 
that $C$ is a {\em precise refinement} of $D$ if $C = \{ U_1,\ldots,U_n \}$ and $D = \{ V_1,\ldots,V_n \}$ with
$U_i \subseteq V_i$ for each $i=1,\ldots,n$.

Here is our main criterion for the existence of star-refinements.  Its statement and proof are
motivated by \cite[\S2.4, Lemma 2.2]{Howes}.
In the statement, we write $V^c$ for the complement of a set $V$.

\begin{lemma}[Shrinking Lemma]
\label{ShrinkingLemma}
Suppose $\C$ is a collection of finite coverings of $X$ satisfying the following two axioms:
\begin{itemize}
\item[(S1)] If $C_1,C_2 \in \C$ then $C_1 \cap C_2 \in \C$.
\item[(S2)] If $D = \{ U_1,\ldots,U_n \} \in \C$, then there exists a refinement
$C = \{ V_1,\ldots,V_n \}$ of $D$ 
with $V_i \subseteq U_i$ for all $i=1,\ldots,n$ 
such that $\{ U_i , V_i^c \} \in \C$ for each $i=1,\ldots,n$.
\end{itemize}
Then $\C$ is a base for a uniform structure on $X$.
\end{lemma}

Note that $C$ in the statement of the lemma is required to be a covering, but it is not required to belong to $\C$.

\begin{proof}
Since axiom (S1) is identical to (C2)(a), it suffices to verify axiom (C2)(b).
Let $D \in \C$, and for each $i$, write $D_i = \{ U_i, V_i^c \} \in \C$.  If $D'$ denotes the intersection
of the binary covers $D_1,\ldots,D_n$, then it suffices by the equivalence of 
(C2)(b) and $({\rm C2})^{\prime}({\rm b})$
to prove that $D'$ delta-refines $D$.
To see that $D'$ delta-refines $D$, let $x \in X$, so that $x \in V_i$ for some $i$.
We claim that $\Star(x,D') \subseteq U_i$.
Indeed, take any $U' \in D'$ such that $x \in U'$.  By definition, we have either $U' \subseteq U_i$ or 
$U' \subseteq V_i^c$, but since $x \in U' \cap V_i$, we cannot have $U' \subseteq V_i^c$;
therefore $U' \subseteq U_i$.  This proves the claim, which implies that
$D'$ is a delta-refinement of $D$ as desired.
\end{proof}

\begin{corollary}
\label{ShrinkingCor}
Let $X$ be a topological space, and suppose that $\B$ is a collection of open sets which are closed under
finite intersections and form a base for the topology of $X$.  Let 
$\C_{\B}$ denote the collection of all finite coverings of $X$ by elements of $\B$, and 
let $\B^c = \{ U^c \; | \; U \in \B \}$.
Suppose that $\C_{\B}$ is non-empty, and that every covering $C \in \C_{\B}$ has a precise refinement by elements of $\B^c$.
Then $\C_{\B}$ is a base for a uniform structure on $X$ which is compatible
with the given topology.
\end{corollary}

\begin{proof}
The fact that $\C_{\B}$ is a base for a uniform structure $\C$ on $X$ follows from
Lemma~\ref{ShrinkingLemma}.  
From the definition of the uniform topology on $(X,\C)$, one sees that 
the uniform topology is coarser than the original topology, and the
two topologies agree if and only if for each $x \in X$ and each neighborhood $U$ of $x$
in the original topology, there exists a covering $C \in \C$ such that
$\Star(x,C) \subseteq U$.  
Without loss of generality we may suppose that $U \in \B$.
Choose any covering $C' = \{ U, U_1,\ldots,U_n \} \in \C_\B$ containing $U$, and let 
$V$ be an element of $\B^c$ contained in $U$.
(The existence of $V$ follows from the existence of a
precise refinement of $C'$ by elements of $\B^c$.)
Finally, let $C = \{ U, U_1 \cap V^c, \ldots, U_n \cap V^c \}$.
Then $C \in \C_{\B}$, and $U$ is the only element of $C$ containing $x$,
so $\Star(x,C) \subseteq U$ as desired.
\end{proof}

The uniform structures considered in Corollary~\ref{ShrinkingCor}, constructed from a special base for
the topology of $X$, are a direct generalization of Rivera-Letelier's construction in \cite{RLUniform}, 
\S3.3--3.4.
(Although Rivera-Letelier uses the entourage formalism for uniform structures, the two constructions
are essentially the same.)
From this point of view, Corollary~\ref{ShrinkingCor} says that the shrinking property for elements of $\C_\B$ 
implies properties $(B_{II})$ and $(B_{III})$ from \cite{RLUniform}.

\section{Uniform structures on rigid analytic and Berkovich spaces}
\label{RigidAnalysisSection}

\subsection{Background from rigid analysis}

In this section, we recall some facts and definitions from rigid analysis which will be used in
the next section; a basic reference for everything we need is \cite{BGR}.  
We require none of the more subtle aspects
of the theory; the main nontrivial fact we use is the Maximum Modulus Principle.

Let $K$ be a field which is complete with respect to a nontrivial 
non-archimedean absolute value, 
and set
\[
\sqrt{|\Kstar|} = \{ \alpha > 0 \; | \; \alpha^n = |z| \text{ for some } z \in K \text{ and some integer } n \geq 1 \}.
\]

For each integer $n \geq 1$, define $T_n$ to be the ring of restricted power series in $n$ variables over $K$.
(``Restricted'' means that the coefficients tend to zero as the degree of the monomial terms increases.)
A basic fact about $T_n$ is that it becomes a Banach algebra over $K$ when equipped with a natural norm called the {\em Gauss norm}. 
An {\em affinoid algebra} over $K$ is a quotient ring of $T_n$ for some $n\geq 1$.
Using the fact that every ideal in $T_n$ is closed, one shows that 
all affinoid algebras are Banach algebras over $K$.
Although there is not in general a canonical choice of norm on an affinoid algebra $\A$, 
there is a canonical equivalence class of norms on $\A$.  In particular, $\A$ has a canonical
topology, independent of the choice of a continuous surjective homomorphism $T_n \onto \A$.

Let $\A$ be an affinoid algebra over $K$, and let $X = {\rm Max}(\A)$ be the 
set of maximal ideals in $\A$.
One thinks of $\A$ as the ``ring of regular functions on $X$'', where a function $f$
is evaluated at a maximal ideal $\m$ by taking the image of $f$ in the quotient field $\A / \m$.
(If $\m$ corresponds to the point $x \in X$, we write $f(x)$ for the image of $f$ in $\A / \m$.)
A basic fact about affinoid algebras is that $\A / \m$ is a finite extension of $K$, and
therefore the absolute value on $K$ has a unique extension to $\A / \m$.
In particular, there is a well-defined map $|f| : \Max(\A) \to \RR_{\geq 0}$ 
which sends $x \in X$ to $|f(x)|$.

In rigid analysis, one equips $X$ with a Grothendieck topology and a sheaf of regular functions, but we will simply think of
$X$ as a topological space, where a base for the topology on $X$ (called the {\em canonical topology}) 
is given by the {\em open Laurent domains} in $X$.
An {\em open Laurent domain} in $X$ is a set of the form 
\[
\{ x \in X \; : \; |f_i(x)| < \alpha_i, \, |g_j(x)| > \beta_j \}
\]
for some $f_1,\ldots,f_s, g_1,\ldots,g_t \in \A$ and $\alpha_1,\ldots,\alpha_s,\beta_1,\ldots,\beta_t \in \sqrt{|\Kstar|}$.
We will refer to the pair $(X,\A)$ (or, by abuse of notation, just to $X$), 
as an {\em affinoid space}.

Similarly, one can define a {\em closed Laurent domain} in $X$ to be a set of the form 
\[
\{ x \in X \; : \; |f_i(x)| \leq \alpha_i, \, |g_j(x)| \geq \beta_j \} \ .
\]

In rigid analysis, one usually considers only what we are calling closed Laurent domains, and they are referred to simply as 
Laurent domains.  
For us, however, open Laurent domains are the more important objects.
If $K$ is algebraically closed, closed Laurent domains are both open and closed in the canonical topology on $X$
(\cite{BGR}, Theorem 7.2.5/3), so one needs to be careful not to abuse the terminology too much!

From the definitions, one sees that an intersection of a finite number of open (resp. closed) Laurent domains is again an open (resp. closed) Laurent domain,
and the complement of an open (resp. closed) Laurent domain is a finite union of closed (resp. open) Laurent domains.

Laurent domains in $X$ are special cases of rational domains in $X$, which themselves are a special kind of affinoid
subdomain of $X$.  We refer the reader to \S7.2 of \cite{BGR} for the definition of rational and affinoid subdomains, which
are extremely important in rigid analysis but do not play an essential part in the present paper.

If $X$ is the closed unit disc $D \subset K$, an open Laurent domain in $X$ is the same thing as an open disc with finitely many closed discs
removed, an object which Rivera-Letelier in \cite{RLUniform} calls an open affinoid domain.

A non-obvious fact which we need from rigid analysis is the following version of the Maximum Modulus Principle
(c.f. Proposition 6.2.1/4 and Lemma 7.3.4/7 of \cite{BGR}):

\begin{theorem}[Maximum Modulus Principle]
Let $X=\Max(\A)$ be an affinoid space.  If $Y \subseteq X$ is a finite union of closed Laurent subdomains of $X$, then
for every $f \in \A$, the function $|f(x)|$ attains both its minimum and maximum values on $Y$.  Moreover, these values 
belong to $\sqrt{|\Kstar|}$.
\end{theorem}

\subsection{The shrinking property for open Laurent covers}
\label{ShrinkingSection}

Let $X = \Max(\A)$ be an affinoid space over the complete non-archimedean field $K$.
Let $\B_L$ be the set of all open Laurent domains in $X$, and let $\C_{\B_L}$ denote
the collection of all finite coverings of $X$ by elements of $\B_L$.  Finally, let $\C_L$ be
the set of all coverings of $X$ having a refinement in $\C_{\B_L}$.  Our goal in this
section is to give a direct proof of the fact that $\C_L$ is a uniform structure on $X$.
More precisely, we will prove the following fact:

\begin{proposition}[Shrinking property for open Laurent covers]
\label{ShrinkingLaurentProp}
Every finite covering $\{ U_1,\ldots,U_n \}$ of $X$ by open Laurent domains has a precise
refinement $\{ V_1,\ldots,V_n \}$ by closed Laurent domains.
\end{proposition}

Although the proposition can be deduced
from Propositions 7.3.4/8 and 7.2.3/6 of \cite{BGR},
for the convenience of the reader we provide a self-contained direct proof.

\begin{proof}
Let $X = U_1 \cup \cdots \cup U_n$ be a covering of $X$ by open Laurent domains.
By a standard argument (see \cite{BGR}, pp. 283-4 of \S7.2.3 ), there are 
functions $f_j^i,g_j^i \in \A$ ($i=1,\ldots,n, \, j=1,\ldots,N$) such that
\[
U_i = \{ x \in X \; | \; |f_j^i(x)| < 1, \, |g_j^i(x)| > 1 \, \forall \, j = 1,\ldots,N \} \ .
\]
It suffices to show that there exists $\gamma \in \sqrt{|\Kstar|}$ with $|\gamma| < 1$
such that $X = \gamma U_1 \cup U_2 \cup \cdots \cup U_n$, where
\[
\gamma U_i = \{ x \in X \; | \; |f_j^i(x)| < |\gamma|, \, |g_j^i(x)| > |\gamma^{-1}| 
\, \forall \, j = 1,\ldots,N \} \ ,
\]
for then one can inductively construct a covering $\{ \overline{\gamma U_1}, \ldots, \overline{\gamma U_n} \}$
of the desired form, where 
\[
\overline{\gamma U_i} = \{ x \in X \; | \; |f_j^i(x)| \leq |\gamma|, \, |g_j^i(x)| \geq |\gamma^{-1}| 
\, \forall \, j = 1,\ldots,N \} \ .
\]

Fix any $\beta \in \sqrt{|\Kstar|}$ with $|\beta| < 1$, and set
\[
Z_{\beta} = (\beta U_1)^c \cap U_2^c \cap \cdots \cap U_n^c \ .
\]
If $Z_{\beta} = \emptyset$, then we're done (take $\gamma = \beta$).
So without loss of generality, we may assume that $Z_{\beta}$ is a non-empty finite union of 
closed Laurent subdomains of $X$.
By the Maximum Modulus Principle, for each $1 \leq j \leq N$ the maximum value of $|f^1_j|$
(resp. the minimum value of $|g^1_j|$) is achieved on $Z_\beta$ and belongs to
$\sqrt{|\Kstar|}$.  Set
\[
\epsilon = \sup_{\substack{ z \in Z_\beta \\ j \in \{ 1,\ldots,N \} } } 
\left( |g^1_j(z)|^{-1},|f^1_j(z)| \right) 
\]
and choose $\alpha \in \sqrt{|\Kstar|}$ with $|\alpha| = \epsilon$.
Since $Z_{\beta} \subseteq (U_2 \cup \cdots \cup U_n)^c \subseteq U_1$, it follows from 
the Maximum Modulus Principle that $0 < \epsilon < 1$.

Since $Z_{\beta} \subseteq (\beta U_1)^c$, we have $|\beta| \leq \epsilon$.
Thus 
\[
0 < |\beta| \leq |\alpha| < 1 \ .
\]
Now choose $\gamma \in \sqrt{|\Kstar|}$ with $|\alpha| < |\gamma| < 1$.  Then
\[
Z_{\gamma} \subseteq Z_{\alpha} \subseteq \overline{\alpha U_1} \ .
\]
But also $Z_{\gamma} \subseteq (\gamma U_1)^c$, so that
\[
Z_{\gamma} \subseteq \overline{\alpha U_1} \cap (\gamma U_1)^c = \emptyset \ .
\]

Therefore $Z_{\gamma} = \emptyset$, which means that 
$X = \gamma U_1 \cup U_2 \cup \cdots \cup U_n$ as desired.
\end{proof}

According to Proposition~\ref{ShrinkingLaurentProp}, every finite covering of $X$ by
elements of $\B_L$ has a precise refinement by elements of $\B_L^c$.
Since the intersection of two open Laurent domains is again an open Laurent domain, 
Proposition~\ref{ShrinkingLaurentProp} and Corollary~\ref{ShrinkingCor} applied to the
collection $\C_{\B_L}$ imply:

\begin{corollary}
$\C_L$ is a uniform structure on $X$ compatible with the canonical topology.
\end{corollary}

It follows immediately from the definitions that the Hausdorff uniform space $(X,\C_L)$ is totally bounded, and
therefore the completion of $(X,\C_L)$ is compact.  
In particular, this shows (independently of Berkovich's theory) that $X$ can be embedded in a natural way as a
dense subspace of a compact Hausdorff topological space $\Xhat$.

\subsection{Comparison with Berkovich's theory}
\label{BerkovichSection}

As in the previous section, let $X = \Max(\A)$ be an affinoid space.
We now explain (in slightly more detail than was done in the introduction) Berkovich's construction 
from \cite{BerkovichBook} of the analytic space $X^{\an}$ associated to $X$.

At the level of points, Berkovich defines $X^{\an}$ to be the set of all 
continuous multiplicative seminorms on $\A$.
A multiplicative seminorm on $\A$ is a function
$| \cdot |_x : \A \rightarrow \RR_{\ge 0}$ such that $|0|_x = 0$,
$|1|_x = 1$, $|f \cdot g|_x = |f|_x \cdot |g|_x$
and $|f + g|_x \le \max \{ |f|_x,|g|_x \}$ for all $f, g \in \A$.  
(Continuity refers to the canonical topology on $\A$.)

The topology on $X^{\an}$ is defined to be the coarsest one for which all functions of the form
$f \mapsto |f|_x$ are continuous.  
A fundamental system of open neighborhoods for this topology is given by the {\em open Laurent domains} in $X^{\an}$ of the form
\[
\{ x \in X^{\an} \; | \; |f_i|_x < \alpha_i, \, |g_j|_x > \beta_j \}
\]
for some $f_1,\ldots,f_s, g_1,\ldots,g_t \in \A$ and $\alpha_1,\ldots,\alpha_s,\beta_1,\ldots,\beta_t \in \sqrt{|\Kstar|}$.
Berkovich proves in \cite{BerkovichBook} that $X^{\an}$ is a
compact Hausdorff space, and it is path-connected if $X$ is connected in the sense of rigid analysis.
Note that since $X^{\an}$ is compact, it comes equipped with a canonical uniform structure.

The space $X = \Max(\A)$ of maximal ideals of $\A$ is naturally embedded in $X^{\an}$ via
the map which sends $x \in X$ to the multiplicative seminorm $f \mapsto |f(x)|$.
This is an embedding of $X$ onto a dense subspace of $X^{\an}$, a fact which one deduces
from the observation that the open Laurent domains in $X$ are precisely the restrictions to $X$ of the open Laurent domains
in $X^{\an}$.

Recall that the space $(\Xhat,\hat{\C}_L)$ discussed in the previous
section is a complete uniform space containing $(X,\C_L)$ as a dense uniform subspace.
Since $X^{\an}$ is also a complete uniform space containing $(X,\C_L)$ as a dense uniform subspace,
it follows from the uniqueness of completions that
$X^{\an}$ and $(\Xhat,\hat{\C}_L)$ are isomorphic as uniform spaces (and in particular also as 
topological spaces).


\medskip

{\bf Remarks:}

\medskip

1. {\em Another point of view.} 
One can reverse the point of view taken in this paper.
Suppose we start with Berkovich's definition of $X^{\an}$, and we take as given that
$X^{\an}$ is a compact Hausdorff space for which the collection $\B^{\an}_L$ of open Laurent domains 
forms a fundamental system
of open neighborhoods.  Then we can deduce directly from general principles (without using 
Proposition~\ref{ShrinkingLaurentProp}) that the collection of coverings $\C_L$ defined in 
\S\ref{ShrinkingSection} is a uniform structure on $X$.  Indeed, as discussed in \S\ref{CompletenessSection},
the canonical uniform structure on $X^{\an}$ is generated by the set of all finite coverings of $X^{\an}$ by elements of 
$\B^{\an}_L$.  Since the restriction of this uniform structure to $X$ is precisely $\C_L$, 
and since the restriction of a uniform structure to a subspace again satisfies axioms (C1) and (C2),
it follows 
that $\C_L$ is a uniform structure on $X$.
But it is still interesting that one can proceed instead as we have done, 
proving {\em directly} that $\C_L$ is a uniform structure on $X$ and deducing from general principles 
that its completion $\hat{X}$ is compact and coincides with the space defined by Berkovich.

\medskip

2. {\em Alternate descriptions of the uniform structure on $X$.} 
By the same reasoning as above, the uniform structure on an affinoid space $X$ can also be described as the collection of all coverings 
having a finite refinement by rational (resp. affinoid) subdomains of $X$.  The reason for using Laurent domains
in \S\ref{ShrinkingSection} is that they are simpler to work with than 
rational or affinoid domains, and they suffice for our purposes.  However, the notion of ``Laurent domain'' is not
transitive, so for globalizing the theory just described it may be better to work with rational or affinoid domains.

\medskip

3. {\em Globalization.}
If $X$ is a quasi-separated rigid analytic space with an admissible finite covering $X = X_1 \cup \cdots \cup X_n$ by affinoids, 
then the Berkovich space $X^{\an}$ associated to $X$ is compact, and the same reasoning as in Remark 1
shows that there is a canonical uniform structure on $X^{\an}$.  By uniqueness of the uniform structure for compact spaces, the uniform structure 
on $X^{\an}$ can be described concretely as the collection of coverings whose restriction to each $X_i$ is a uniform covering.  
So the construction of $X^{\an}$
as the completion of a canonical uniform structure on $X$ works with almost no modification in this 
greater generality.
For example, this gives an alternate way to construct the Berkovich space associated to an arbitrary proper scheme 
over $K$. 

It would be interesting to know if
there is a canonical uniform structure on $X$ whose completion yields
the Berkovich analytification $X^{\an}$ 
for any quasi-separated rigid analytic space $X$ over $K$ having an admissible affinoid covering
of finite type (c.f. \S1.6 of \cite{BerkovichIHES}).

\medskip

4. {\em Relation to the work of Schneider and Van der Put.}
It is worthwhile to compare the uniform space point of view with related earlier work of Schneider and van der Put
\cite{SvdP}.  In \cite{SvdP}, the authors prove that given an affinoid space $X = {\rm Max}(\A)$ over $K$, the space of {\em maximal filters} on the 
underlying $G$-topology of $X$ coincides with the Berkovich analytic space $X^{\an}$.  
By contrast, we propose viewing $X$ as the space of {\em minimal Cauchy filters} on the set $X$ itself.
Although these two points of view are related, it is not immediately clear how to formulate a dictionary between them.  
It would be nice to have a concrete correspondence between the two points of view.  
While the construction from \cite{SvdP} is certainly natural from the point of view of rigid analysis, it relies
on the non-standard notion of a filter on a $G$-topology, whereas our construction involves only the more classical notion of a Cauchy filter on a uniform space.
In addition, working with actual filters on sets, rather than filters on $G$-topologies, 
feels more in the spirit of Berkovich's theory.


\appendix
\section{Definition of uniform structures via entourages}

The definition of a uniform structure in terms of entourages (taken from Bourbaki \cite{BourbakiGT}) is as follows.

\begin{definition}
A {\em uniform structure} $\U$ on a set $X$ is a collection of subsets $V$ of $X \times X$, called {\em entourages}, which form
a filter and satisfy the following additional axioms:
\begin{itemize}
\item[(E1)] Each $V \in \U$ contains the diagonal in $X \times X$.
\item[(E2)] If $V \in \U$ then $V' := \{ (x,y) \in X \times X \; | \; (y,x) \in V \} \in \U$.
\item[(E3)] If $V \in \U$, then there exists $W \in \U$ such that
\[
W \circ W := \{ (x,y) \in X \times X \; | \; \exists \, z \in X \text{ with } (x,z) \in W \text{ and } (z,y) \in W \}
\]
is contained in $V$.
\end{itemize}
\end{definition}

A {\em uniform space} is a set $X$ together with a uniform structure on $X$.  

Every metric space has a canonical uniform structure generated by the entourages $V_\epsilon = \{ x,y \in X \; | \; d(x,y) < \epsilon \}$ for $\epsilon > 0$,
and every uniform space is a topological space by declaring that the open neighborhoods of a point $x \in X$ are the sets of the form
$U[x] := \{ y \in X \; | (x,y) \in U \}$ for $U \in \U$. 

The definitions of a uniform structure in terms of uniform covers or entourages are equivalent.
The connection between the two definitions is as follows.
Given a uniform space $(X,\U)$ in the entourage sense, define a covering $C$ of $X$ to be a uniform cover if
there exists an entourage $U \in \U$ such that for each $x \in X$, there is a $V \in C$ with
$U[x] \subseteq V$.
Conversely, given a uniform space $(X,\C)$ in the uniform cover sense, define the entourages to be the supersets of
$\bigcup_{V \in C} V \times V$ as $C$ ranges over the uniform covers of $X$.
These two procedures are inverse to one another,
and furnish an equivalence between the two notions of uniform space.

\bibliographystyle{plain}
\bibliography{uniform}

\end{document}